# The asymptotic structure of nearly unstable non-negative integer-valued AR(1) models

FEIKE C. DROST[1,*], RAMON VAN DEN AKKER[2] and BAS J.M. WERKER[1,**]

[1]*Econometrics & Finance group, CentER, Tilburg University, The Netherlands, P.O. Box 90153, 5000 LE Tilburg, The Netherlands.*
*E-mail:* [*]*F.C.Drost@TilburgUniversity.nl;* [**]*B.J.M.Werker@TilburgUniversity.nl*
[2]*Econometrics group, CentER, Tilburg University, The Netherlands, P.O. Box 90153, 5000 LE Tilburg, The Netherlands. E-mail: R.vdnAkker@TilburgUniversity.nl*

This paper considers non-negative integer-valued autoregressive processes where the autoregression parameter is close to unity. We consider the asymptotics of this 'near unit root' situation. The local asymptotic structure of the likelihood ratios of the model is obtained, showing that the limit experiment is Poissonian. To illustrate the statistical consequences we discuss efficient estimation of the autoregression parameter and efficient testing for a unit root.

## 1. Introduction

The non-negative integer-valued autoregressive process of the order 1 (INAR(1)) was introduced by Al-Osh and Alzaid (1987) and Alzaid (1988) as a non-negative integer-valued analogue of the AR(1) process. Al-Osh and Alzaid (1990) and Du and Li (1991) extended this work to INAR($p$) processes. Recently there has been a growing interest in INAR processes. Without going into details we mention some recent (theoretical) contributions to the literature on INAR processes: Freeland and McCabe (2005), Jung, Ronning and Tremayne (2005), Silva and Oliveira (2005), Silva and Silva (2006), Zheng, Basawa and Datta (2006), Neal and Subba Rao (2007) and Drost, Van den Akker and Werker (2008a, 2008b). Applications of INAR processes in the medical sciences can be found in, for example, Franke and Seligmann (1993), Bélisle *et al.* (1998) and Cardinal, Roy and Lambert (1999); an application to psychometrics in Böckenholt (1999a), an application to environmentology in Thyregod et al. (1999); recent applications to economics in, for example, Böckenholt (1999b), Berglund and Brännäs (2001), Brännäs and Hellström (2001), Rudholm (2001), Böckenholt (2003), Brännäs and Quoreshi (2004), Freeland and McCabe (2004), Gourieroux and Jasiak (2004) and Drost, Van den

Akker and Werker (2008c); and Ahn, Gyemin and Jongwoo (2000) and Pickands III and Stine (1997) considered queueing applications.

This paper considers a nearly nonstationary INAR(1) model and derives its limit experiment (in the Le Cam framework). Our main result is that this limit experiment is Poissonian. This is surprising since limit experiments are usually locally asymptotically quadratic (LAQ; see Jeganathan (1995), Le Cam and Yang (1990) and Ling and McAleer (2003)) and even non-regular models often enjoy a shift structure (see Hirano and Porter (2003a, 2003b)), whereas the Poisson limit experiment does not enjoy these two properties. The result is indeed surprising since Ispány, Pap and van Zuijlen (2003a) established a functional limit theorem with an Ornstein–Uhlenbeck limit process from which one would conjecture a standard LAQ-type limit experiment. On a technical level the proof of the convergence to a Poisson limit experiment is interesting, since the 'score' can be split into two parts that have different rates of convergence. To illustrate the statistical consequences of the convergence to a Poisson limit experiment, we exploit this limit experiment to construct efficient estimators of the autoregression parameter and to construct an efficient test for the null hypothesis of a unit root. Since the INAR(1) process is a particular branching process with immigration, this also partially solves the question (see Wei and Winnicki (1990)) of how to estimate the offspring mean efficiently. Furthermore, we show that the ordinary least squares (OLS) estimator, considered by Ispány, Pap and van Zuijlen (2003a, 2003b, 2005), is inefficient. Surprisingly, the OLS estimator even has a lower rate of convergence than the efficient one. Related to this, we show that the classical Dickey–Fuller test for a unit root has no power against local alternatives induced by the limit experiment. More precisely, as we will see below, the autoregressive parameter in these local alternatives is of the form $1 - h/n^2$ with $h > 0$ and $n$ denoting the number of observations. Of course, for alternatives at a further distance the Dickey–Fuller test will have power but the efficient test can perfectly discriminate between the null and the alternative in such a case.

An INAR(1) process (starting at 0) is defined by the recursion, $X_0 = 0$, and,

$$X_t = \vartheta \circ X_{t-1} + \varepsilon_t, \qquad t \in \mathbb{N}, \tag{1}$$

where (by definition an empty sum equals 0),

$$\vartheta \circ X_{t-1} = \sum_{j=1}^{X_{t-1}} Z_j^{(t)}. \tag{2}$$

Here $(Z_j^{(t)})_{j \in \mathbb{N}, t \in \mathbb{N}}$ is a collection of i.i.d. Bernoulli variables with success probability $\theta \in [0, 1]$, independent of the i.i.d. innovation sequence $(\varepsilon_t)_{t \in \mathbb{N}}$ with distribution $G$ on $\mathbb{Z}_+ = \mathbb{N} \cup \{0\}$. All these variables are defined on a probability space $(\Omega, \mathcal{F}, \mathbb{P}_{\theta, G})$. If we work with fixed $G$, we usually drop the subscript $G$. The random variable $\vartheta \circ X_{t-1}$ is called the binomial thinning of $X_{t-1}$ (this operator was introduced by Steutel and van Harn (1979)) and, conditionally on $X_{t-1}$, it follows a binomial distribution with success probability $\theta$ and a number of trials equal to $X_{t-1}$. Equation (1) can be interpreted as a branching process with immigration. The outcome $X_t$ is composed of $\vartheta \circ X_{t-1}$, the elements of $X_{t-1}$ that survive during $(t-1,t]$, and $\varepsilon_t$, the number of immigrants during $(t-1,t]$. Here the number of immigrants is independent of the survival of elements of $X_{t-1}$. Moreover, each element of $X_{t-1}$ survives with probability $\theta$ and its survival has no effect on the survival of the other elements. From a statistical point of view, the difference between the literature on INAR processes and the literature on branching processes with immigration is that in the latter setting one commonly observes both the $X$ process and the $\varepsilon$ process, whereas one only observes the $X$ process in the INAR setting, which complicates inference drastically. Compared to the familiar AR(1) processes, inference for INAR(1) processes is also more complicated, since, even if $\theta$ is known, observing $X$ does not imply observing $\varepsilon$. From the definition of an INAR process it immediately follows that $\mathbb{E}_{\theta,G}[X_t|X_{t-1},\ldots,X_0] = \theta X_{t-1} + \mathbb{E}_G \varepsilon_1$, which (partially) explains the 'AR' in 'INAR'. It is well known that, if $\theta \in [0,1)$ and $\mathbb{E}_G \varepsilon_1 < \infty$, which is called the 'stable' case, there exists an initial distribution, $\nu_{\theta,G}$, such that $X$ is stationary if $\mathcal{L}(X_0) = \nu_{\theta,G}$. Of course, the INAR(1) process is non-stationary if $\theta = 1$: under $\mathbb{P}_1$ the process $X$ is nothing but a standard random walk with drift on $\mathbb{Z}_+$ (but note that $X$ is nondecreasing under $\mathbb{P}_1$). We call this situation 'unstable' or say that the process has a 'unit root'. Although the unit root is on the boundary of the parameter space, it is an important parameter value since in many applications the estimates of $\theta$ are close to 1.

Denote the law of $(X_0,\ldots,X_n)$ under $\mathbb{P}_{\theta,G}$ on the measurable space $(\mathcal{X}_n, \mathcal{A}_n) = (\mathbb{Z}_+^{n+1}, 2^{\mathbb{Z}_+^{n+1}})$ by $\mathbb{P}_{\theta,G}^{(n)}$. For $G$ known, the global model of interest is thus $(\mathbb{P}_{\theta,G} \mid \theta \in [0,1])$. The model restricted to the stable case $\theta \in [0,1)$, has been shown to be locally asymptotically normal (LAN) in Drost, Van den Akker and Werker (2008b) and Section 4.3.2 in Van den Akker (2007). For this stable case, the OLS estimator is consistent and asymptotically normal. The focus of interest of the present paper is, however, the unstable case $\theta = 1$. Therefore we will introduce the local parameter $h \geq 0$ and take the autoregressive parameter $\theta_n = 1 - h/n^2$ in (2). This is formalized below.

In our applications we mainly consider two sets of assumptions on $G$: (i) $G$ is known or (ii) $G$ is completely unknown (apart from some regularity conditions). For expository reasons, let us, for the moment, focus on the case that $G$ is completely known and the goal is to estimate $\theta$. We use 'local-to-unity' asymptotics to take the 'increasing statistical difficulty' in the neighborhood of the unit root into account, that is, we consider local alternatives to the unit root in such a way that the increasing degree of difficulty to discriminate between these alternatives and the unit root compensates the increase of information contained in the sample as the number of observations grows. This approach is well known; it originated from the work of Chan and Wei (1987) and Philips (1987), who studied the behavior of a given estimator (OLS) in a nearly unstable AR(1) setting, and Jeganathan (1995), whose results yield the asymptotic structure of nearly unstable AR models. Following this approach, we introduce the sequence of nearly unstable INAR experiments $\mathcal{E}_n(G) = (\mathcal{X}_n, \mathcal{A}_n, (\mathbb{P}_{1-h/n^2}^{(n)} \mid h \geq 0))$, $n \in \mathbb{N}$. The 'localizing rate' $n^2$ will become apparent later on. It is surprising that the localizing rate is $n^2$, since for the classical nearly unstable AR(1) model one has rate $n\sqrt{n}$ (non-zero intercept) or $n$ (no intercept). Suppose that we have found an estimator $\hat{h}_n$ with 'nice properties'; then this corresponds to the estimate $\hat{\theta}_n = 1 - \hat{h}_n/n^2$ of $\theta$ in the global experiment of interest.

300 *F.C. Drost, R. van den Akker and B.J.M. Werker*

To our knowledge, Ispány, Pap and van Zuijlen (2003a) were the first to study estimation in a nearly unstable INAR(1) model. These authors study the behavior of the OLS estimator and they use a localizing rate $n$ instead of $n^2$. However, $n^2$ is the proper localizing rate and, in Proposition 4.3, we show indeed that the OLS estimator is an exploding estimator in $(\mathcal{E}_n(G))_{n\in\mathbb{N}}$, that is, it has not even the 'right' rate of convergence. The question then arises how we should estimate $h$. Instead of analyzing the asymptotic behavior of a given estimator, we derive the asymptotic structure of the experiments themselves by determining the limit experiment (in the Le Cam sense) of $(\mathcal{E}_n(G))_{n\in\mathbb{N}}$. This limit experiment gives bounds to the accuracy of inference procedures and suggests how to construct efficient ones.

The main contribution of this paper is to determine the limit experiment of $(\mathcal{E}_n(G))_{n\in\mathbb{N}}$. Remember that (see, e.g., Le Cam (1986), Le Cam and Yang (1990), Van der Vaart (1991), Shiryaev and Spokoiny (1999) or Van der Vaart (2000) Chapter 9), the sequence of experiments $(\mathcal{E}_n(G))_{n\in\mathbb{N}}$ is said to converge to a limit experiment (in Le Cam's weak topology) $\mathcal{E} = (\mathcal{X}, \mathcal{A}, (\mathbb{Q}_h \mid h \geq 0))$ if, for every finite subset $I \subset \mathbb{R}_+$ and every $h_0 \in \mathbb{R}_+$, we have

$$\left(\frac{\mathrm{d}\mathbb{P}^{(n)}_{1-h/n^2}}{\mathrm{d}\mathbb{P}^{(n)}_{1-h_0/n^2}}(X_0,\ldots,X_n)\right)_{h\in I} \xrightarrow{d} \left(\frac{\mathrm{d}\mathbb{Q}_h}{\mathrm{d}\mathbb{Q}_{h_0}}(Z)\right)_{h\in I}, \qquad \text{under } \mathbb{P}^{(n)}_{1-h_0/n^2}.$$

To see that it is indeed reasonable to expect $n^2$ as the proper localizing rate we briefly discuss the case of geometrically distributed innovations (in the remainder we treat general $G$). In case $G = \text{Geometric}(1/2)$, that is, $G$ puts mass $(1/2)^{k+1}$ at $k \in \mathbb{Z}_+$, it is an easy exercise to verify for $h > 0$ (the geometric distribution allows us, using Newton's binomial formula, to obtain explicit expressions for the transition probabilities from $X_{t-1}$ to $X_t$ if $X_t \geq X_{t-1}$),

$$\frac{\mathrm{d}\mathbb{P}^{(n)}_{1-h,r_n}}{\mathrm{d}\mathbb{P}^{(n)}_1}(X_0,\ldots,X_0)$$

$$\xrightarrow{p} \begin{cases} 0, & \text{if } \dfrac{r_n}{n^2} \to 0, \\ \exp\left(-\dfrac{hG(0)\mathbb{E}_G\varepsilon_1}{2}\right) = \exp\left(-\dfrac{h}{4}\right), & \text{if } \dfrac{r_n}{n^2} \to 1, \\ 1, & \text{if } \dfrac{r_n}{n^2} \to \infty, \end{cases} \qquad \text{under } \mathbb{P}_1.$$

This has two important implications: First, it indicates that $n^2$ is indeed the proper localizing rate. Intuitively, if we go faster than $n^2$ we cannot distinguish $\mathbb{P}^{(n)}_{1-h/r_n}$ from $\mathbb{P}^{(n)}_1$, and if we go slower we can distinguish $\mathbb{P}^{(n)}_{1-h/r_n}$ perfectly from $\mathbb{P}^{(n)}_1$. Second, since $\exp(-h/4) < 1$ we have, by Le Cam's first lemma, no contiguity of $\mathbb{P}^{(n)}_{1-h/n^2}$ with respect to $\mathbb{P}^{(n)}_1$. (Remark 2 after Theorem 2.1 gives an example of sets that yield this non-contiguity.) This lack of contiguity is unfortunate for several reasons. Most important, if we had contiguity the limiting behavior of $(d\mathbb{P}^{(n)}_{1-h/n^2}/d\mathbb{P}^{(n)}_1)_{h\in I}$ would determine the limit experiment, whereas we now need to consider the behavior of $(d\mathbb{P}^{(n)}_{1-h/n^2}/d\mathbb{P}^{(n)}_{1-h_0/n^2})_{h\in I}$ for all $h_0 \geq 0$. And it implies that the global sequence of experiments does not have the common LAQ structure (see Jeganathan (1995)) at $\theta = 1$. This differs from the traditional AR(1) process $Y_0 = 0$, $Y_t = \mu + \theta Y_{t-1} + u_t$, $u_t$ i.i.d. $N(0, \sigma^2)$, with $\mu \neq 0$ and $\sigma^2$ known, that enjoys this LAQ property at $\theta = 1$: the limit experiment at $\theta = 1$ is the usual normal location experiment (i.e., the model is LAN) and the localizing rate is $n^{3/2}$. The limit experiment at $\theta = 1$ for $Y_0 = 0$, $Y_t = \theta Y_{t-1} + u_t$, $u_t$ i.i.d. $N(0, \sigma^2)$, with $\sigma^2$ known, does not have the LAN structure; the limit experiment is of the locally asymptotically Brownian functional type (a special class of LAQ experiments; see Jeganathan (1995)) and the localizing rate is $n$. Thus although the INAR(1) process and the traditional AR(1) process both are walks with drift at $\theta = 1$, their statistical properties 'near $\theta = 1$' are very different. In Section 3 we prove that the limit experiment of $(\mathcal{E}_n(G))_{n\in\mathbb{N}}$ corresponds to one draw from a Poisson distribution with mean $hG(0)\mathbb{E}_G\varepsilon_1/2$, that is,

$$\mathcal{E}(G) = \left(\mathbb{Z}_+, 2^{\mathbb{Z}_+}, \left(\text{Poisson}\left(\frac{hG(0)\mathbb{E}_G\varepsilon_1}{2}\right) \,\Big|\, h \geq 0\right)\right).$$

We indeed recognize $\exp(-hG(0)\mathbb{E}_G\varepsilon_1/2)$ as the likelihood ratio at $h$ relative to $h_0 = 0$ in the experiment $\mathcal{E}(G)$. Due to the lack of enough smoothness of the likelihood ratios around the unit root, this convergence of experiments is not obtained by the usual (generally applicable) techniques, but rather by a direct approach. Since the transition probability is the convolution of a binomial distribution with $G$ and the fact that certain binomial experiments converge to a Poisson limit experiment, one might be tempted to think that the convergence $\mathcal{E}_n(G) \to \mathcal{E}(G)$ follows, in some way, from this convergence. As is clear from the proof of Theorem 3.1 this reasoning is not valid.

The remainder of the paper is organized as follows: In Section 2 we discuss some preliminary properties that provide insight into the behavior of a nearly unstable INAR(1) process and are needed in the rest of the paper. The main result is stated and proved in Section 3. Section 4 uses our main result to analyze some estimation and testing problems. We consider efficient inference of $h$, the deviation from a unit root, in the nearly unstable case for two settings. The first setting, discussed in Section 4.1, treats the case that the immigration distribution $G$ is completely known. The second setting, analyzed in Section 4.2, considers a semi-parametric model, where hardly any conditions on $G$ are imposed. Furthermore, we show in Section 4.1 that the OLS estimator is explosive under the local alternative $\theta_n = 1 - h/n^2$. Finally, we discuss testing for a unit root in Section 4.3. We show that the traditional Dickey–Fuller test has no (local) power, but that an intuitively obvious test is efficient. Appendix A contains some auxiliary results and Appendix B gathers some proofs.

## 2. Preliminaries

This section discusses some basic properties of nearly unstable INAR(1) processes. Besides giving insight into the behavior of a nearly unstable INAR(1) process, these proper-ties are a key input in the next sections. To enhance readability the proofs are organized in Appendix B.

First, we introduce the following notation: The mean of $\varepsilon_t$ is denoted by $\mu_G$ and its variance by $\sigma_G^2$. The probability mass function corresponding to $G$, the distribution of the innovations $\varepsilon_t$, is denoted by $g$. The probability mass function of the binomial distribution with parameters $\theta \in [0,1]$ and $n \in \mathbb{Z}_+$ is denoted by $b_{n,\theta}$.

Given $X_{t-1} = x_{t-1}$, the random variables $\varepsilon_t$ and $\vartheta \circ X_{t-1}$ are independent and $X_{t-1} - \vartheta \circ X_{t-1}$, 'the number of deaths during $(t-1,t]$', follows a binomial$(X_{t-1}, 1-\theta)$ distribution. This interpretation yields the following representation of the transition probabilities,

$$\begin{aligned} P^\theta_{x_{t-1},x_t} &= \mathbb{P}_\theta\{X_t = x_t | X_{t-1} = x_{t-1}\} \\ &= \sum_{k=0}^{x_{t-1}} \mathbb{P}_\theta\{\varepsilon_t = x_t - x_{t-1} + k, X_{t-1} - \vartheta \circ X_{t-1} = k | X_{t-1} = x_{t-1}\} \\ &= \sum_{k=0}^{x_{t-1}} b_{x_{t-1}, 1-\theta}(k) g(\Delta x_t + k), \end{aligned}$$

where $\Delta x_t = x_t - x_{t-1}$, and $g(i) = 0$ for $i < 0$. Under $\mathbb{P}_1$ we have $X_t = \mu_G t + \sum_{i=1}^t (\varepsilon_i - \mu_G)$, and $P^1_{x_{t-1},x_t} = g(\Delta x_t)$, $x_{t-1}, x_t \in \mathbb{Z}_+$. Hence, under $\mathbb{P}_1$, an INAR(1) process is nothing but an integer-valued random walk with drift.

The next proposition is basic, but often applied in the sequel.

**Proposition 2.1.** *If $\sigma_G^2 < \infty$, we have for $h \geq 0$,*

$$\lim_{n \to \infty} \mathbb{E}_{1-h/n^2} \left[ \frac{1}{n^2} \sum_{t=1}^n X_t - \frac{\mu_G}{2} \right]^2 = 0. \tag{3}$$

*If $\sigma_G^2 < \infty$, then we have for $\alpha > 0$ and every sequence $(\theta_n)_{n \in \mathbb{N}}$ in $[0,1]$,*

$$\lim_{n \to \infty} \frac{1}{n^{3+\alpha}} \sum_{t=1}^n \mathbb{E}_{\theta_n} X_t^2 = 0. \tag{4}$$

*Remark 1.* Convergence in probability for the case $h > 0$ in (3) cannot be concluded from the convergence in probability in (3) for $h = 0$ by contiguity arguments. The reason is (see Remark 2 after the proof of Theorem 2.1) that $\mathbb{P}^{(n)}_{1-h/n^2}$ is not contiguous with respect to $\mathbb{P}^{(n)}_1$.

Next, we consider the thinning process $(\vartheta \circ X_{t-1})_{t \geq 1}$. Under $\mathbb{P}_{1-h/n^2}$, $X_{t-1} - \vartheta \circ X_{t-1}$, conditional on $X_{t-1}$, a binomial$(X_{t-1}, h/n^2)$ distribution follows. So we expect that, near unity, many 'deaths' do not occur in any time interval $(t-1,t]$. The following proposition gives a precise statement where we use the following notation: For $h \geq 0$ and $n \in \mathbb{N}$,

$$A_n^h = \left\{ z \in \mathbb{Z}_+ \,\Big|\, \frac{h(z+1)}{n^2} < \frac{1}{2} \right\}, \qquad \mathcal{A}_n^h = \{(X_0, \ldots, X_{n-1}) \in A_n^h \times \cdots \times A_n^h\}. \tag{5}$$

These sets are introduced for the following reasons: By Proposition A.1 we have for $x \in A_n^h$ $\sum_{k=r}^x \mathrm{b}_{x,h/n^2}(k) \leq 2\,\mathrm{b}_{x,h/n^2}(r)$ for $r = 2, 3$ and terms of the form $(1 - \frac{h}{n^2})^{-2}$ can be bounded neatly without having to make statements of the form 'for $n$ large enough', or having to refer to 'up to a constant depending on $h$'. Also, recall the notation $\Delta X_t = X_t - X_{t-1}$.

**Proposition 2.2.** *If $\sigma_G^2 < \infty$, then we have for all sequences $(\theta_n)_{n \in \mathbb{N}}$ in $[0,1]$ and for all $h \geq 0$,*

$$\lim_{n \to \infty} \mathbb{P}_{\theta_n}(\mathcal{A}_n^h) = 1. \tag{6}$$

*Moreover, if $\sigma_G^2 < \infty$ and $h \geq 0$, we have,*

$$\lim_{n \to \infty} \mathbb{P}_{1-h/n^2}\{\exists t \in \{1, \ldots, n\} : X_{t-1} - \vartheta \circ X_{t-1} \geq 2\} = 0. \tag{7}$$

Finally, we derive the limit distribution of the number of downward movements of $X$ during $[0, n]$. The probability that the Bernoulli variable $1\{\Delta X_t < 0\}$ equals one is small. Intuitively, the dependence over time of this indicator process is not too strong, so it is not unreasonable to expect that a 'Poisson law of small numbers' holds. As the following theorem shows, this is indeed the case.

**Theorem 2.1.** *If $\sigma_G^2 < \infty$, then we have for $h \geq 0$,*

$$\sum_{t=1}^n 1\{\Delta X_t < 0\} \xrightarrow{d} \mathrm{Poisson}\left(\frac{hg(0)\mu_G}{2}\right), \qquad under\ \mathbb{P}_{1-h/n^2}. \tag{8}$$

**Remark 2.** Since $\sum_{t=1}^n 1\{\Delta X_t < 0\}$ equals zero under $\mathbb{P}_1^{(n)}$ and converges in distribution to a non-degenerated limit under $\mathbb{P}_{1-h/n^2}^{(n)}$ ($h > 0$, $0 < g(0) < 1$), we see that $\mathbb{P}_{1-h/n^2}^{(n)}$ is not contiguous with respect to $\mathbb{P}_1^{(n)}$ for $h > 0$.

## 3. The limit experiment: one observation from a Poisson distribution

For easy reference, we introduce the following assumption.

**Assumption 3.1.** *A probability distribution $G$ on $\mathbb{Z}_+$ is said to satisfy Assumption 3.1 if one of the following two conditions holds:*

(1) $\mathrm{Support}(G) = \{0,\ldots,M\}$ for some $M \in \mathbb{N}$;
(2) $\mathrm{Support}(G) = \mathbb{Z}_+$, $\sigma_G^2 < \infty$ and $g$ is eventually decreasing, that is, there exists $M \in \mathbb{N}$ such that $g(k+1) \leq g(k)$ for $k \geq M$.

The rest of this section is devoted to the following theorem.

**Theorem 3.1.** *Suppose $G$ satisfies Assumption 3.1. Then the limit experiment of $(\mathcal{E}_n(G))_{n \in \mathbb{N}}$ is given by*

$$\mathcal{E}(G) = (\mathbb{Z}_+, 2^{\mathbb{Z}_+}, (\mathbb{Q}_h \mid h \geq 0)),$$

*with $\mathbb{Q}_h = \mathrm{Poisson}(hg(0)\mu_G/2)$. More precisely, for $h \geq 0$ and $h_0 > 0$ we have, under $\mathbb{P}^{(n)}_{1-h_0/n^2}$,*

$$\frac{\mathrm{d}\mathbb{P}^{(n)}_{1-h/n^2}}{\mathrm{d}\mathbb{P}^{(n)}_{1-h_0/n^2}}(X_0,\ldots,X_n) \xrightarrow{d} \frac{\mathrm{d}\mathbb{Q}_h}{\mathrm{d}\mathbb{Q}_{h_0}}(Z) = \exp\left(-\frac{(h-h_0)g(0)\mu_G}{2}\right)\left(\frac{h}{h_0}\right)^Z, \qquad (9)$$

*while for $h \geq 0$ and $h_0 = 0$ we have, under $\mathbb{P}^{(n)}_1$,*

$$\frac{\mathrm{d}\mathbb{P}^{(n)}_{1-h/n^2}}{\mathrm{d}\mathbb{P}^{(n)}_1}(X_0,\ldots,X_n) \xrightarrow{d} \frac{\mathrm{d}\mathbb{Q}_h}{\mathrm{d}\mathbb{Q}_0}(Z) = \exp\left(-\frac{hg(0)\mu_G}{2}\right)\mathbf{1}\{Z=0\}. \qquad (10)$$

**Proof.** Introduce for $h,h_0 \geq 0$,

$$\mathcal{L}_n(h,h_0) = \log\frac{\mathrm{d}\mathbb{P}^{(n)}_{1-h/n^2}}{\mathrm{d}\mathbb{P}^{(n)}_{1-h_0/n^2}} = \sum_{t=1}^n \log\frac{P^{1-h/n^2}_{X_{t-1},X_t}}{P^{1-h_0/n^2}_{X_{t-1},X_t}}.$$

Note, if $\sum_{t=1}^n \mathbf{1}\{\Delta X_t < 0\} > 0$ and $h_0 > 0$, that $\mathcal{L}_n(0,h_0) = -\infty$ and thus $\mathrm{d}\mathbb{P}^{(n)}_0/\mathrm{d}\mathbb{P}^{(n)}_{1-h_0/n^2} = 0$. Because $\mathcal{L}_n(h,h_0)$ is complicated to analyze, split the transition probability $P^{1-h/n^2}_{x_{t-1},x_t}$ into a leading term,

$$L_n(x_{t-1},x_t,h) = \begin{cases} \displaystyle\sum_{k=-\Delta x_t}^{-\Delta x_t+1} b_{x_{t-1},h/n^2}(k)g(\Delta x_t+k), & \text{if } \Delta x_t < 0, \\ \displaystyle\sum_{k=0}^{1} b_{x_{t-1},h/n^2}(k)g(\Delta x_t+k), & \text{if } \Delta x_t \geq 0, \end{cases}$$

and a remainder term,

$$R_n(x_{t-1},x_t,h) = \begin{cases} \displaystyle\sum_{k=-\Delta x_t+2}^{x_{t-1}} b_{x_{t-1},h/n^2}(k)g(\Delta x_t+k), & \text{if } \Delta x_t < 0, \\ \displaystyle\sum_{k=2}^{x_{t-1}} b_{x_{t-1},h/n^2}(k)g(\Delta x_t+k), & \text{if } \Delta x_t \geq 0. \end{cases}$$

We introduce a simpler version of $\mathcal{L}_n(h, h_0)$ in which the remainder terms are removed,

$$\tilde{\mathcal{L}}_n(h, h_0) = \sum_{t=1}^n \log \frac{L_n(X_{t-1}, X_t, h)}{L_n(X_{t-1}, X_t, h_0)}.$$

The difference between $\tilde{\mathcal{L}}_n(h, h_0)$ and $\mathcal{L}_n(h, h_0)$ is negligible.

**Lemma 3.1.** *If $G$ satisfies Assumption 3.1, then we have for $h \geq 0$ and $h_0 \geq 0$,*

$$\tilde{\mathcal{L}}_n(h, h_0) = \mathcal{L}_n(h, h_0) + \mathrm{o}(\mathbb{P}_{1-h_0/n^2}; 1). \tag{11}$$

To enhance readability the proof of the lemma is organized in Appendix B. Hence, the limit distribution of the random vector $(\mathcal{L}_n(h, h_0))_{h \in I}$, for a finite subset $I \subset \mathbb{R}_+$, is the same as the limit distribution of $(\tilde{\mathcal{L}}_n(h, h_0))_{h \in I}$. It easily follows, using (7), that $\tilde{\mathcal{L}}_n(h, h_0)$ can be decomposed as

$$\tilde{\mathcal{L}}_n(h, h_0) = \sum_{t=1}^n \frac{X_{t-1} - 2}{n^2} \log\left(\frac{1 - h/n^2}{1 - h_0/n^2}\right)^{n^2} + S_n^+(h, h_0)$$
$$+ S_n^-(h, h_0) + \mathrm{o}(\mathbb{P}_{1-h_0/n^2}; 1), \tag{12}$$

where $S_n^+(h, h_0) = \sum_{t: \Delta X_t \geq 0} W_{tn}^+$ and $S_n^-(h, h_0) = \sum_{t: \Delta X_t = -1} W_{tn}^-$ are defined by (here $\sum_{t: \Delta X_t = -1}$ is shorthand for $\sum_{1 \leq t \leq n: \Delta X_t = -1}$, and for $\sum_{t: \Delta X_t \geq 0}$ the same convention is used),

$$W_{tn}^+ = \log\left[\frac{g(\Delta X_t)(1 - h/n^2)^2 + X_{t-1}(h/n^2)(1 - h/n^2)g(\Delta X_t + 1)}{g(\Delta X_t)(1 - h_0/n^2)^2 + X_{t-1}(h_0/n^2)(1 - h_0/n^2)g(\Delta X_t + 1)}\right],$$
$$W_{tn}^- = \log\left[\frac{X_{t-1}(h/n^2)(1 - h/n^2)g(0) + (X_{t-1}(X_{t-1} - 1)/2)(h^2/n^4)g(1)}{X_{t-1}(h_0/n^2)(1 - h_0/n^2)g(0) + (X_{t-1}(X_{t-1} - 1)/2)(h_0^2/n^4)g(1)}\right].$$

So we need to determine the asymptotic behavior of the terms in (12). By (3) we have,

$$\log\left[\left(\frac{1 - h/n^2}{1 - h_0/n^2}\right)^{n^2}\right] \frac{1}{n^2} \sum_{t=1}^n (X_{t-1} - 2) \xrightarrow{p} -\frac{(h - h_0)\mu_G}{2}, \qquad \text{under } \mathbb{P}_{1-h_0/n^2}. \tag{13}$$

The next lemma yields the behavior of $S_n^+(h, h_0)$, the second term of (12); see Appendix B for the proof.

**Lemma 3.2.** *If $G$ satisfies Assumption 3.1, then we have for $h \geq 0$ and $h_0 \geq 0$,*

$$S_n^+(h, h_0) \xrightarrow{p} \frac{(h - h_0)(1 - g(0))\mu_G}{2}, \qquad \text{under } \mathbb{P}_{1-h_0/n^2}. \tag{14}$$

Finally, we discuss the term $S_n^-(h, h_0)$ in (12). Under $\mathbb{P}_1$ this term is not present, so we only need to consider $h_0 > 0$. We organize the result in the following lemma; see Appendix B for the proof.

**Lemma 3.3.** *If $G$ satisfies Assumption 3.1, then we have for $h_0 > 0$ and $h \geq 0$,*

$$S_n^-(h, h_0) = \log\left[\frac{h}{h_0}\right] \sum_{t=1}^n 1\{\Delta X_t < 0\} + o(\mathbb{P}_{1-h_0/n^2}; 1), \tag{15}$$

*where we set $\log(0) = -\infty$ and $\log(0) \cdot 0 = 0$.*

To complete the proof of the theorem, note that we obtain from Lemmas 3.1–3.3, (12) and (13):

$$\mathcal{L}_n(h, h_0) = \tilde{\mathcal{L}}_n(h, h_0) + o(\mathbb{P}_{1-h_0/n^2}; 1)$$
$$= -\frac{(h-h_0)g(0)\mu_G}{2} + \log\left[\frac{h}{h_0}\right] \sum_{t=1}^n 1\{\Delta X_t < 0\} + o(\mathbb{P}_{1-h_0/n^2}; 1),$$

where we interpret $\log(0) = -\infty$, $\log(0) \cdot 0 = 0$ and $\log(h/0) \sum_{t=1}^n 1\{\Delta X_t < 0\} = 0$ when $h_0 = 0$, $h > 0$. Hence, Theorem 2.1 implies that, for a finite subset $I \subset \mathbb{R}_+$,

$$(\mathcal{L}_n(h, h_0))_{h \in I} \xrightarrow{d} \left(\log \frac{d\mathbb{Q}_h}{d\mathbb{Q}_{h_0}}(Z)\right)_{h \in I}, \qquad \text{under } \mathbb{P}_{1-h_0/n^2},$$

which concludes the proof. □

**Remark 3.** In the proof we have seen that,

$$\log \frac{d\mathbb{P}_{1-h/n^2}^{(n)}}{d\mathbb{P}_{1-h_0/n^2}^{(n)}} = -\frac{(h-h_0)g(0)\mu_G}{2} + \log\left[\frac{h}{h_0}\right] \sum_{t=1}^n 1\{\Delta X_t < 0\} + o(\mathbb{P}_{1-h_0/n^2}; 1).$$

So, heuristically, we can interpret $\sum_{t=1}^n 1\{\Delta X_t < 0\}$ as an 'approximately sufficient statistic'.

## 4. Applications

This section considers the following applications as an illustration of the statistical consequences of the convergence of experiments. We discuss the efficient estimation of $h$, the deviation from a unit root, in the nearly unstable case for two settings. The first setting, discussed in Section 4.1, treats the case that $G$ is completely known. And the second setting, analyzed in Section 4.2, considers a semi-parametric model, where hardly any conditions on $G$ are imposed. Finally, we discuss testing for a unit root in Section 4.3.

## 4.1. Efficient estimation of *h* in nearly unstable INAR models (*G* known)

In this section $G$ is assumed to be known. So we consider the sequence of experiments $(\mathcal{E}_n(G))_{n\in\mathbb{N}}$. As before, we denote the observation from the limit experiment $\mathcal{E}(G)$ by $Z$, and $\mathbb{Q}_h = \mathrm{Poisson}(hg(0)\mu_G/2)$.

Since we have established convergence of $(\mathcal{E}_n(G))_{n\in\mathbb{N}}$ to $\mathcal{E}(G)$, an application of the Le Cam–Van der Vaart asymptotic representation theorem yields the following proposition.

**Proposition 4.1.** *Suppose $G$ satisfies Assumption 3.1. If $(T_n)_{n\in\mathbb{N}}$ is a sequence of estimators of $h$ in the sequence of experiments $(\mathcal{E}_n(G))_{n\in\mathbb{N}}$ such that $\mathcal{L}(T_n|\mathbb{P}_{1-h/n^2}) \to Z_h$ for all $h \geq 0$, then there exists a map $t:\mathbb{Z}_+ \times [0,1] \to \mathbb{R}$ such that $Z_h = \mathcal{L}(t(Z,U)|\mathbb{Q}_h \times \mathrm{Uniform}[0,1])$ (i.e., $U$ is distributed uniformly on $[0,1]$ and independent of the observation $Z$ from the limit experiment $\mathcal{E}(G)$).*

**Proof.** The sequence $(\mathcal{E}_n(G))_{n\in\mathbb{N}}$ converges to the experiment $\mathcal{E}(G)$ (Theorem 3.1). Since $\mathcal{E}(G)$ is dominated by counting measure on $\mathbb{Z}_+$, the result follows by applying the Le Cam–Van der Vaart asymptotic representation theorem (see, e.g., Van der Vaart (1991), Theorem 3.1, or Van der Vaart (2000), Theorem 9.3). □

Thus, for any set of limit laws of an estimator there is a randomized estimator in the limit experiment that has the same set of laws. If the asymptotic performance of an estimator is considered to be determined by its sets of limit laws, the limit experiment thus gives a lower bound to what is possible: Along the sequence of experiments you cannot do better than the best procedure in the limit experiment. To discuss efficient estimation we need to prescribe what we judge to be optimal in the Poisson limit experiment. Often a normal location experiment is the limit experiment. For such a normal location experiment, that is, estimate $h$ on the basis of one observation $Y$ from $\mathrm{N}(h,\tau)$ ($\tau$ known), it is natural to restrict to location-equivariant estimators. For this class one has a convolution property (see, e.g., Bickel et al. (1998), Van der Vaart (2000) or Wong (1992)): the law of every location-equivariant estimator $T$ of $h$ can be decomposed as $T \stackrel{d}{=} Y + V$, where $V$ is independent of $Y$. This yields, by Anderson's lemma (see, e.g., Lemma 8.5 in Van der Vaart (2000)), efficiency of $Y$ (within the class of location-equivariant estimators) for all bowl-shaped loss functions. To be more general, there are convolution results for shift experiments. However, the Poisson limit experiment $\mathcal{E}(G)$ does not have a natural shift structure. In such a Poisson setting it seems reasonable to minimize variance amongst the unbiased estimators. See Ling and McAleer (2003) for a similar approach for LAQ limit experiments.

**Definition 4.1.** *An estimator $\hat{h}$ for $h$ is called efficient in the experiment $\mathcal{E}(G)$ if it is unbiased, that is, $\mathbb{E}_h \hat{h} = h$ for all $h \geq 0$, and minimizes the variance amongst all unbiased (randomized) estimators of $h$.*

The next proposition is an immediate consequence of the Lehmann–Scheffé theorem.

**Proposition 4.2.** *If $0 < g(0) < 1$ and $\mu_G < \infty$, then $2Z/g(0)\mu_G$ is an efficient estimator of $h$ in the experiment $\mathcal{E}(G)$. The variance of this estimator, under $\mathbb{Q}_h$, is given by $2h/g(0)\mu_G$.*

A combination with Proposition 4.1 yields a variance lower bound to asymptotically unbiased estimators in the sequence of experiments $(\mathcal{E}_n(G))_{n\in\mathbb{N}}$.

**Corollary 4.1.** *Suppose $G$ satisfies Assumption 3.1. If $(T_n)_{n\in\mathbb{N}}$ is an estimator of $h$ in the sequence of experiments $(\mathcal{E}_n(G))_{n\in\mathbb{N}}$ such that $\mathcal{L}(T_n|\mathbb{P}_{1-h/n^2}) \to Z_h$ with $\int z\,\mathrm{d}Z_h(z) = h$ for all $h \geq 0$, then we have*

$$\int (z-h)^2 \,\mathrm{d}Z_h(z) \geq \frac{2h}{g(0)\mu_G} \qquad \text{for all } h \geq 0. \tag{16}$$

It is not unnatural to restrict to estimators that satisfy $\mathcal{L}(T_n|\mathbb{P}_{1-h/n^2}) \to Z_h$. We make the additional restriction that $\int z\,\mathrm{d}Z_h(Z) = h$, that is, the *limit* distribution is unbiased. Now, based on the previous proposition, it is natural to call an estimator in this class efficient if it attains the variance bound (16). To demonstrate the efficiency of a given estimator, one only needs to show that it belongs to the class of asymptotically unbiased estimators, and that it attains the bound. How should we estimate $h$? Recall, that we interpreted $\sum_{t=1}^n 1\{\Delta X_t < 0\}$ as an approximately sufficient statistic for $h$. Hence, it is natural to try to construct an efficient estimator based on this statistic. Using Theorem 2.1 we see that this is indeed possible.

**Corollary 4.2.** *If Assumption 3.1 holds, then the estimator,*

$$\hat{h}_n = \frac{2\sum_{t=1}^n 1\{\Delta X_t < 0\}}{g(0)\mu_G}, \tag{17}$$

*is an efficient estimator of $h$ in the sequence $(\mathcal{E}_n(G))_{n\in\mathbb{N}}$.*

Finally, we discuss the commonly used OLS estimator when $\theta_n = 1 - h/n^2$. Rewriting $X_t = \vartheta \circ X_{t-1} + \varepsilon_t = \mu_G + \theta_n X_{t-1} + u_t$ for $u_t = \varepsilon_t - \mu_G + \vartheta \circ X_{t-1} - \theta_n X_{t-1}$, we obtain the autoregression $X_t - \mu_G = \theta_n X_{t-1} + u_t$, which can also be written as $n^2(X_t - X_{t-1} - \mu_G) = h(-X_{t-1}) + n^2 u_t$ (note that indeed $\mathbb{E}_{\theta_n} u_t = \mathbb{E}_{\theta_n} X_{t-1} u_t = 0$). So the OLS estimator of $\theta_n$ is given by

$$\hat{\theta}_n^{\mathrm{OLS}} = \frac{\sum_{t=1}^n X_{t-1}(X_t - \mu_G)}{\sum_{t=1}^n X_{t-1}^2}, \tag{18}$$

and the OLS estimator of $h$ is given by

$$\hat{h}_n^{\mathrm{OLS}} = -\frac{n^2 \sum_{t=1}^n X_{t-1}(X_t - X_{t-1} - \mu_G)}{\sum_{t=1}^n X_{t-1}^2} = n^2(1 - \hat{\theta}_n^{\mathrm{OLS}}).$$

Ispány, Pap and van Zuijlen (2003a) showed that $n^{3/2}(\hat{\theta}_n^{\text{OLS}} - \gamma_n) \xrightarrow{d} N(0, \sigma_\gamma^2)$ under $\mathbb{P}_{\gamma_n}$ for $\gamma_n = 1 - h_n/n$, $h_n \to \gamma$, and $\sigma_\gamma^2$ depending on $\gamma$. This means that the OLS estimator can be used to distinguish alternatives at rate $n$. Since the convergence of experiments takes place at rate $n^2$, the OLS estimator deteriorates under the localizing rate $n^2$.

**Proposition 4.3.** *If $\mathbb{E}_G \varepsilon_1^3 < \infty$, then we have for all $h \geq 0$,*

$$|\hat{h}_n^{\text{OLS}}| \xrightarrow{p} \infty, \qquad \text{under } \mathbb{P}_{1-h/n^2}.$$

*Remark 4.* A similar result holds for the OLS estimator in the model where $G$ is unknown.

Thus the OLS estimator cannot distinguish local alternatives at rate $n^2$; at lower rates (up to $n^{3/2}$) it is capable of distinguishing alternatives. In this sense it does not have the right rate of convergence.

## 4.2. Efficient estimation of *h* in nearly unstable INAR models (*G* unknown)

So far we have assumed that $G$ is known. In this section, where we instead consider a semi-parametric model, we hardly impose conditions on $G$ (see, e.g., Bickel and Kwon (2001) or Wefelmeyer (1996) for general theories on semi-parametric stationary Markov models, and Drost, Klaassen and Werker (1997) for group-based time series models). The dependence of $\mathbb{P}_\theta$ upon $G$ is made explicit by adding a subscript: $\mathbb{P}_{\theta,G}$. Formally, we consider the sequence of experiments,

$$\mathcal{E}_n = (\mathbb{Z}_+^{n+1}, 2^{\mathbb{Z}_+^{n+1}}, (\mathbb{P}_{1-h/n^2,G}^{(n)} \mid (h,G) \in [0,\infty) \times \mathcal{G})), \qquad n \in \mathbb{N},$$

where $\mathcal{G}$ is the set of all distributions on $\mathbb{Z}_+$ that satisfy Assumption 3.1.

The goal is to estimate $h$ efficiently. Here efficient, just as in the previous section, means asymptotically unbiased with minimal variance. Since the semi-parametric model is more realistic, the estimation of $h$ becomes more difficult. As we will see, the situation for our semi-parametric model is quite fortunate: we can estimate $h$ with the same asymptotic precision as in the case where $G$ is known. In semi-parametric statistics this is called adaptive estimation.

The efficient estimator for the case where $G$ is known cannot be used anymore, since it depends on $g(0)$ and $\mu_G$. The obvious idea is to replace these objects by estimates. The next proposition provides consistent estimators.

**Proposition 4.4.** *Let $h \geq 0$ and $G$ satisfy $\sigma_G^2 < \infty$. Then we have*

$$\hat{g}_n(0) = \frac{1}{n} \sum_{t=1}^n 1\{X_t = X_{t-1}\} \xrightarrow{p} g(0) \quad \text{and} \quad \hat{\mu}_{G,n} = \frac{X_n}{n} \xrightarrow{p} \mu_G, \qquad \text{under } \mathbb{P}_{1-h/n^2,G}.$$

**Proof.** Notice first that we have

$$\frac{1}{n}\sum_{t=1}^{n}(X_{t-1} - \vartheta \circ X_{t-1}) \xrightarrow{p} 0, \qquad \text{under } \mathbb{P}_{1-h/n^2,G}, \tag{19}$$

thus condition on $X_{t-1}$ and use (3),

$$0 \leq \frac{1}{n}\sum_{t=1}^{n}\mathbb{E}_{1-h/n^2,G}(X_{t-1} - \vartheta \circ X_{t-1}) = \frac{h}{n^3}\sum_{t=1}^{n}\mathbb{E}_{1-h/n^2,G}X_{t-1} \to 0.$$

Using that $|1\{X_t = X_{t-1}\} - 1\{\varepsilon_t = 0\}| = 1$ only if $X_{t-1} - \vartheta \circ X_{t-1} \geq 1$, we easily obtain by using (19),

$$\left|\hat{g}_n(0) - \frac{1}{n}\sum_{t=1}^{n}1\{\varepsilon_t = 0\}\right| \leq \frac{1}{n}\sum_{t=1}^{n}1\{X_{t-1} - \vartheta \circ X_{t-1} \geq 1\} \leq \frac{1}{n}\sum_{t=1}^{n}(X_{t-1} - \vartheta \circ X_{t-1}) \xrightarrow{p} 0.$$

Now the result for $\hat{g}_n(0)$ follows by applying the weak law of large numbers to $n^{-1}\sum_{t=1}^{n}1\{\varepsilon_t = 0\}$. Next, consider $\hat{\mu}_{G,n}$. We have, using (19) and the weak law of large numbers for $n^{-1}\sum_{t=1}^{n}\varepsilon_t$,

$$\hat{\mu}_{G,n} = \frac{X_n}{n} = \frac{1}{n}\sum_{t=1}^{n}(X_t - X_{t-1})$$

$$= \frac{1}{n}\sum_{t=1}^{n}\varepsilon_t - \frac{1}{n}\sum_{t=1}^{n}(X_{t-1} - \vartheta \circ X_{t-1}) \xrightarrow{p} \mu_G, \qquad \text{under } \mathbb{P}_{1-h/n^2,G},$$

which concludes the proof. □

From the previous proposition we have $\hat{h}_n - \tilde{h}_n \xrightarrow{p} 0$, under $\mathbb{P}_{1-h/n^2,G}$, where

$$\tilde{h}_n = \frac{2\sum_{t=1}^{n}1\{\Delta X_t < 0\}}{\hat{g}_n(0)\hat{\mu}_{G,n}}.$$

This implies that estimation of $h$ in the semi-parametric experiments $(\mathcal{E}_n)_{n\in\mathbb{N}}$ is not harder than the estimation of $h$ in $(\mathcal{E}_n(G))_{n\in\mathbb{N}}$. In semi-parametric parlor: The semi-parametric problem is adaptive to $\mathcal{G}$. The precise statement is given in the following corollary; the proof is trivial.

**Corollary 4.3.** *If $(T_n)_{n\in\mathbb{N}}$ is a sequence of estimators in the semi-parametric sequence of experiments $(\mathcal{E}_n)_{n\in\mathbb{N}}$ such that $\mathcal{L}(T_n|\mathbb{P}_{1-h/n^2,G}) \to Z_{h,G}$ with $\int z\,dZ_{h,G}(z) = h$ for all $(h,G) \in [0,\infty) \times \mathcal{G}$, then we have*

$$\int (z-h)^2\,dZ_{h,G}(z) \geq \frac{2h}{g(0)\mu_G} \qquad \textit{for all } (h,G) \in [0,\infty) \times \mathcal{G}.$$

*The estimator $\tilde{h}_n$ satisfies the conditions and achieves the variance bound.*

### 4.3. Testing for a unit root

This section discusses testing for a unit root in an INAR(1) model. We consider the case where $G$ is known and satisfies Assumption 3.1. We want to test the hypothesis $H_0 : \theta = 1$ versus $H_1 : \theta < 1$. Hellstrom (2001) considered this problem from the perspective that one wants to use standard (i.e., OLS) software routines. He derives, by Monte Carlo simulations, the finite sample null distributions for a Dickey–Fuller test of a random walk with Poisson distributed errors. This (standard) Dickey–Fuller test statistic is given by the *usual* (i.e., non-corrected) $t$-test that the slope parameter equals 1, that is,

$$\tau_n = \frac{\hat{\theta}_n^{\mathrm{OLS}} - 1}{\sqrt{\sigma_G^2 (\sum_{t=1}^n X_{t-1}^2)^{-1}}},$$

where $\hat{\theta}_n^{\mathrm{OLS}}$ is given by (18). Under $H_0$, that is, under $\mathbb{P}_1$, we have $\tau_n \xrightarrow{d} \mathrm{N}(0,1)$. To analyze the power of this test, and since $\mathcal{E}_n(G) \to \mathcal{E}(G)$, we consider the performance of $\tau_n$ along the sequence $\mathcal{E}_n(G)$. The following proposition shows, however, that the asymptotic probability that the null hypothesis is rejected remains $\alpha$ for all alternatives. Obviously, this does not exclude power of the Dickey–Fuller test under local alternatives at rate $n^{3/2}$ (which indeed it has).

**Proposition 4.5.** *If $\mathbb{E}_G \varepsilon_1^3 < \infty$, we have for all $h \geq 0$,*

$$\tau_n \xrightarrow{d} \mathrm{N}(0,1), \qquad \text{under } \mathbb{P}_{1-h/n^2}, \qquad \text{which yields } \lim_{n \to \infty} \mathbb{P}_{1-h/n^2}(\text{reject } H_0) = \alpha.$$

**Proof.** From Ispány, Pap and van Zuijlen (2003a) the result easily follows. □

We propose the intuitively obvious tests

$$\psi_n(X_0, \ldots, X_n) = \begin{cases} \alpha, & \text{if } \sum_{t=1}^n 1\{\Delta X_t < 0\} = 0, \\ 1, & \text{if } \sum_{t=1}^n 1\{\Delta X_t < 0\} \geq 1, \end{cases}$$

that is, reject $H_0$ if the process ever moves down and reject $H_0$ with probability $\alpha$ if there are no downward movements. We will see that this obvious test is, in fact, efficient. To discuss the efficiency of tests, we recall the implication of the Le Cam–Van der Vaart asymptotic representation theorem to testing (see Theorem 7.2 in Van der Vaart (1991)). Let $\alpha \in (0,1)$ and $\phi_n$ be a sequence of tests in $(\mathcal{E}_n(G))_{n \in \mathbb{N}}$ such that $\limsup_{n \to \infty} \mathbb{E}_1 \phi_n(X_0, \ldots, X_n) \leq \alpha$. Then we have

$$\limsup_{n \to \infty} \mathbb{E}_{1-h/n^2} \phi_n(X_0, \ldots, X_n) \leq \sup_{\phi \in \Phi_\alpha} \mathbb{E}_h \phi(Z) \qquad \text{for all } h > 0,$$

where $\Phi_\alpha$ is the collection of all level $\alpha$ tests for testing $H_0 : h = 0$ versus $H_1 : h > 0$ in the Poisson limit experiment $\mathcal{E}(G)$. If we have equality in the previous display, it is natural to call a test $\phi_n$ efficient. It is obvious that the uniform most powerful test in the Poisson limit experiment is given by

$$\phi(Z) = \begin{cases} \alpha, & \text{if } Z = 0, \\ 1, & \text{if } Z \geq 1. \end{cases}$$

Its power function is given by $\mathbb{E}_0 \phi(Z) = \alpha$ and $\mathbb{E}_h \phi(Z) = 1 - (1-\alpha)\exp(-hg(0)\mu_G/2)$. Using Theorem 2.1 we find

$$\lim_{n \to \infty} \mathbb{E}_1 \psi_n(X_0, \ldots, X_n) = \alpha$$

and

$$\lim_{n \to \infty} \mathbb{E}_{1-h/n^2} \psi_n(X_0, \ldots, X_n) = 1 - (1-\alpha)\exp\left(-\frac{hg(0)\mu_G}{2}\right) \quad \text{for } h > 0.$$

We conclude that the test $\psi_n$ is indeed efficient.

## Appendix A: Auxiliaries

The following result is basic (see, e.g., Feller (1968), pages 150–151), but since it is heavily applied, we state it here for convenience.

**Proposition A.1.** *Let $m \in \mathbb{N}$, $p \in (0,1)$. If $r > mp$, we have*

$$\sum_{k=r}^{m} b_{m,p}(k) \leq b_{m,p}(r) \frac{r(1-p)}{r - mp}. \tag{20}$$

*So, if $1 > mp$, we have for $r = 2, 3$,*

$$\sum_{k=r}^{m} b_{m,p}(k) \leq 2 b_{m,p}(r). \tag{21}$$

For convenience we recall Theorem 1 in Serfling (1975).

**Lemma A.1.** *Let $Z_1, \ldots, Z_n$ (possibly dependent) 0–1 valued random variables and set $S_n = \sum_{t=1}^n Z_t$. Let $Y$ be Poisson distributed with mean $\sum_{t=1}^n \mathbb{E} Z_t$. Then we have*

$$\sup_{A \subset \mathbb{Z}_+} |\mathbb{P}\{S_n \in A\} - \mathbb{P}\{Y \in A\}| \leq \sum_{t=1}^n (\mathbb{E} Z_t)^2 + \sum_{t=1}^n \mathbb{E}|\mathbb{E}[Z_t | Z_1, \ldots, Z_{t-1}] - \mathbb{E} Z_t|.$$

## Appendix B: Proofs

**Proof of Proposition 2.1.** Obviously $\mathrm{Var}_1(\sum_{t=1}^n X_t) = \mathrm{O}(n^3)$ and $\lim_{n\to\infty} n^{-2} \times \sum_{t=1}^n \mathbb{E}_1 X_t = \mu_G/2$, which yields (3) for $h=0$. Next, we prove (3) for $h>0$. Straightforward calculations show, for $\theta < 1$,

$$\mathbb{E}_\theta \sum_{t=1}^n X_t = \mu_G \sum_{t=1}^n \frac{1-\theta^t}{1-\theta} = \mu_G \left[ \frac{n}{1-\theta} - \frac{\theta - \theta^{n+1}}{(1-\theta)^2} \right],$$

which yields

$$\lim_{n\to\infty} \frac{1}{n^2} \mathbb{E}_{1-h/n^2} \sum_{t=1}^n X_t$$
$$= \lim_{n\to\infty} \frac{\mu_G}{n^2} \left[ \frac{n}{h/n^2} - \frac{1 - h/n^2 - [1-(n+1)h/n^2 + ((n+1)n/2)h^2/n^4 + \mathrm{o}(1/n^2)]}{h^2/n^4} \right]$$
$$= \frac{\mu_G}{2}. \tag{1}$$

To treat the variance of $n^{-2}\sum_{t=1}^n X_t$, we use the following simple relations; see also Ispány, Pap and van Zuijlen (2003a), for $0 < \theta < 1$, $s, t \geq 1$,

$$\mathrm{Cov}_\theta(X_t, X_s) = \theta^{|t-s|} \mathrm{Var}_\theta X_{t\wedge s},$$
$$\mathrm{Var}_\theta X_t = \frac{1-\theta^{2t}}{1-\theta^2}\sigma_G^2 + \frac{(\theta-\theta^t)(1-\theta^t)}{1-\theta^2}\mu_G \leq (\sigma_G^2 + \mu_G)\frac{1-\theta^{2t}}{1-\theta^2}. \tag{2}$$

From this we obtain, as $n \to \infty$,

$$\mathrm{Var}_{1-h/n^2}\left(\frac{1}{n^2}\sum_{t=1}^n X_t\right) = \frac{1}{n^4}\sum_{t=1}^n \left(1 + 2\sum_{s=t+1}^n \left(1-\frac{h}{n^2}\right)^{s-t}\right) \mathrm{Var}_{1-h/n^2} X_t$$
$$\leq \frac{1}{n}2n(\sigma_G^2+\mu_G)\frac{1}{n^2}\frac{1}{1-(1-h/n^2)^2}\frac{1}{n}\sum_{t=1}^n \left(1-\left(1-\frac{h}{n^2}\right)^{2t}\right) \to 0.$$

Together with (1) this completes the proof of (3) for $h > 0$. To prove (4), note that $X_t \leq \sum_{i=1}^t \varepsilon_i$. Hence $\mathbb{E}_{\theta_n} X_t^2 \leq \mathbb{E}_1 X_t^2 = \sigma_G^2 t + \mu_G^2 t^2$, which yields the desired conclusion. □

**Proof of Proposition 2.2.** Equation (6) easily follows since, for a sequence $(\theta_n)_{n\in\mathbb{N}}$ in $[0,1]$, (4) implies

$$\mathbb{P}_{\theta_n}\{\exists\, 0 \leq t \leq n : X_t > n^{7/4}\} \leq \frac{1}{n^{7/2}}\sum_{t=1}^n \mathbb{E}_{\theta_n} X_t^2 \to 0 \quad \text{as } n \to \infty. \tag{3}$$

To obtain (7) note that, for $X_{t-1} \in A_n^h$ we have, using the bound (21),

$$\mathbb{P}_{1-h/n^2}\{X_{t-1} - \vartheta \circ X_{t-1} \geq 2 | X_{t-1}\} = \sum_{k=2}^{X_{t-1}} b_{X_{t-1},h/n^2}(k) \leq 2\,b_{X_{t-1},h/n^2}(2) \leq \frac{h^2 X_{t-1}^2}{n^4}.$$

By (4) this yields,

$$\lim_{n\to\infty} \mathbb{P}_{1-h/n^2}(\{\exists t \in \{1,\ldots,n\} : X_{t-1} - \vartheta \circ X_{t-1} \geq 2\} \cap \mathcal{A}_n^h) \leq \lim_{n\to\infty} \frac{h^2}{n^4} \sum_{t=1}^n \mathbb{E}_{1-h/n^2} X_{t-1}^2 = 0.$$

Since we already showed $\lim_{n\to\infty} \mathbb{P}_{1-h/n^2}(\mathcal{A}_n^h) = 1$, this yields (7). $\square$

**Proof of Theorem 2.1.** If $g(0) = 0$, then $\Delta X_t < 0$ implies $X_{t-1} - \vartheta \circ X_{t-1} \geq 2$. Hence, (7) implies $\sum_{t=1}^n 1\{\Delta X_t < 0\} \xrightarrow{p} 0$ under $\mathbb{P}_{1-h/n^2}$. Since the Poisson distribution with mean 0 concentrates all its mass at 0, this yields the result. The cases $h = 0$ or $g(0) = 1$ (recall $X_0 = 0$) are also trivial. So we consider the case $h > 0$ and $0 < g(0) < 1$. For notational convenience, abbreviate $\mathbb{P}_{1-h/n^2}$ by $\mathbb{P}_n$ and $\mathbb{E}_{1-h/n^2}$ by $\mathbb{E}_n$. Put $Z_t = 1\{\Delta X_t = -1, \varepsilon_t = 0\}$ and notice that $0 \leq 1\{\Delta X_t < 0\} - Z_t = 1\{\Delta X_t \leq -2\} + 1\{\Delta X_t = -1, \varepsilon_t \geq 1\}$. From (7) it now follows that

$$0 \leq \sum_{t=1}^n 1\{\Delta X_t < 0\} - \sum_{t=1}^n Z_t \leq 2 \sum_{t=1}^n 1\{X_{t-1} - \vartheta \circ X_{t-1} \geq 2\} \xrightarrow{p} 0, \qquad \text{under } \mathbb{P}_n.$$

Thus it suffices to prove that $\sum_{t=1}^n Z_t \xrightarrow{d} \text{Poisson}(hg(0)\mu_G/2)$ under $\mathbb{P}_n$. We do this by applying Lemma A.1. Introduce random variables $Y_n$, where $Y_n$ follows a Poisson distribution with mean $\lambda_n = \sum_{t=1}^n \mathbb{E}_n Z_t$. And let $Z$ follow a Poisson distribution with mean $hg(0)\mu_G/2$. From Lemma A.1 we obtain the bound

$$\sup_{A \subset \mathbb{Z}_+} \left| \mathbb{P}_n\left\{ \sum_{t=1}^n Z_t \in A \right\} - \Pr\{Y_n \in A\} \right|$$
$$\leq \sum_{t=1}^n (\mathbb{E}_n Z_t)^2 + \sum_{t=1}^n \mathbb{E}_n |\mathbb{E}_n[Z_t - \mathbb{E}_n Z_t | Z_1, \ldots, Z_{t-1}]|.$$

If we prove that

$$\text{(i) } \sum_{t=1}^n (\mathbb{E}_n Z_t)^2 \to 0, \qquad \text{(ii) } \sum_{t=1}^n \mathbb{E}_n Z_t \to \frac{hg(0)\mu_G}{2},$$

$$\text{(iii) } \sum_{t=1}^n \mathbb{E}_n |\mathbb{E}_n[Z_t - \mathbb{E}_n Z_t | Z_1, \ldots, Z_{t-1}]| \to 0,$$

all hold as $n \to \infty$, then the result follows since we then have for all $z \in \mathbb{R}$,

$$\left|\mathbb{P}_n\left\{\sum_{t=1}^n Z_t \leq z\right\} - \Pr(Z \leq z)\right| \leq \left|\mathbb{P}_n\left\{\sum_{t=1}^n Z_t \leq z\right\} - \Pr\{Y_n \leq z\}\right|$$
$$+ |\Pr\{Y_n \leq z\} - \Pr(Z \leq z)| \to 0.$$

First we tackle (i). Using that, conditional on $X_{t-1}$, $\varepsilon_t$ and $X_{t-1} - \vartheta \circ X_{t-1} \sim \text{Bin}_{X_{t-1}, h/n^2}$ being independent, we obtain

$$\mathbb{E}_n Z_t = \mathbb{P}_n\{\varepsilon_t = 0, X_{t-1} - \vartheta \circ X_{t-1} = 1\} = \frac{hg(0)}{n^2} \mathbb{E}_n X_{t-1} \left(1 - \frac{h}{n^2}\right)^{X_{t-1}-1} \leq \frac{hg(0)}{n^2} \mathbb{E}_n X_{t-1}.$$

Then (i) is easily obtained using (4),

$$\lim_{n \to \infty} \sum_{t=1}^n (\mathbb{E}_n Z_t)^2 \leq \lim_{n \to \infty} \frac{h^2 g^2(0)}{n^4} \sum_{t=1}^n \mathbb{E}_n X_{t-1}^2 = 0.$$

Next we consider (ii). If we prove the relation,

$$\lim_{n \to \infty} \left|\frac{1}{n^2} \sum_{t=1}^n \mathbb{E}_n X_{t-1} - \frac{1}{n^2} \sum_{t=1}^n \mathbb{E}_n X_{t-1} \left(1 - \frac{h}{n^2}\right)^{X_{t-1}-1}\right| = 0,$$

it is immediate that (ii) follows from (3). To prove the previous display, we introduce $B_n = \{\forall t \in \{1, \ldots, n\} : X_t \leq n^{7/4}\}$ with $\lim_{n \to \infty} \mathbb{P}_n(B_n) = 1$ (see (3)). On the event $B_n$ we have $n^{-2} X_t \leq n^{-1/4}$ for $t = 1, \ldots, n$. This yields

$$0 \leq \mathbb{E}_n X_{t-1} \left(1 - \left(1 - \frac{h}{n^2}\right)^{X_{t-1}-1}\right) \leq \mathbb{E}_n X_{t-1} \left(1 - \left(1 - \frac{h}{n^2}\right)^{X_{t-1}}\right) 1_{B_n} + \mathbb{E}_n X_{t-1} 1_{B_n^c}$$

$$\leq \mathbb{E}_n \left[1_{B_n} X_{t-1} \sum_{j=1}^{X_{t-1}} \binom{X_{t-1}}{j} \left(\frac{h}{n^2}\right)^j\right] + \mathbb{E}_n X_{t-1} 1_{B_n^c} \leq \frac{1}{n^{1/4}} \exp(h) \mathbb{E}_n X_{t-1} + \mathbb{E}_n X_{t-1} 1_{B_n^c}.$$

Using $\mathbb{P}_n(B_n) \to 1$ and (3) we obtain,

$$\lim_{n \to \infty} \frac{1}{n^2} \sum_{t=1}^n \mathbb{E}_n X_{t-1} 1_{B_n^c} \leq \lim_{n \to \infty} \sqrt{\mathbb{E}_n \left(\frac{1}{n^2} \sum_{t=1}^n X_{t-1}\right)^2 \mathbb{P}_n(B_n^c)} = \sqrt{\left(\frac{\mu_G}{2}\right)^2 \cdot 0} = 0.$$

By (4) we have $\lim_{n \to \infty} n^{-9/4} \sum_{t=1}^n \mathbb{E}_n X_{t-1} = 0$. Combination with the previous two displays yields the result.

Finally, we prove (iii). Let $\mathcal{F}^\varepsilon = (\mathcal{F}_t^\varepsilon)_{t \geq 1}$ and $\mathcal{F}^X = (\mathcal{F}_t^X)_{t \geq 0}$ be the filtrations generated by $(\varepsilon_t)_{t \geq 1}$ and $(X_t)_{t \geq 0}$, respectively, that is, $\mathcal{F}_t^\varepsilon = \sigma(\varepsilon_1, \ldots, \varepsilon_t)$ and $\mathcal{F}_t^X =$

$\sigma(X_0, \ldots, X_t)$. Note that we have, for $t \geq 2$,

$$\mathbb{E}_n |\mathbb{E}_n[Z_t - \mathbb{E}_n Z_t | Z_1, \ldots, Z_{t-1}]|$$
$$\leq \mathbb{E}_n |\mathbb{E}_n[Z_t - \mathbb{E}_n Z_t | \mathcal{F}_{t-1}^\varepsilon, \mathcal{F}_{t-1}^X]| \qquad (4)$$
$$= \frac{hg(0)}{n^2} \mathbb{E}_n \left| X_{t-1} \left(1 - \frac{h}{n^2}\right)^{X_{t-1}-1} - \mathbb{E}_n X_{t-1} \left(1 - \frac{h}{n^2}\right)^{X_{t-1}-1} \right|.$$

Using the reverse triangle inequality we obtain

$$\left| \mathbb{E}_n \left| X_{t-1} \left(1 - \frac{h}{n^2}\right)^{X_{t-1}-1} - \mathbb{E}_n X_{t-1} \left(1 - \frac{h}{n^2}\right)^{X_{t-1}-1} \right| - \mathbb{E}_n |X_{t-1} - \mathbb{E}_n X_{t-1}| \right|$$
$$\leq \mathbb{E}_n \left| X_{t-1} \left(1 - \left(1 - \frac{h}{n^2}\right)^{X_{t-1}-1}\right) - \mathbb{E}_n X_{t-1} \left(1 - \left(1 - \frac{h}{n^2}\right)^{X_{t-1}-1}\right) \right|$$
$$\leq 2 \mathbb{E}_n X_{t-1} \left(1 - \left(1 - \frac{h}{n^2}\right)^{X_{t-1}-1}\right).$$

We have already seen in the proof of (ii) that

$$\lim_{n\to\infty} \frac{1}{n^2} \sum_{t=1}^n \mathbb{E}_n X_{t-1} \left(1 - \left(1 - \frac{h}{n^2}\right)^{X_{t-1}-1}\right) = 0.$$

A combination of the previous two displays with (4) now easily yields the bound

$$\sum_{t=1}^n \mathbb{E}_n |\mathbb{E}_n[Z_t - \mathbb{E}_n Z_t | Z_1, \ldots, Z_{t-1}]| \leq \mathrm{o}(1) + \frac{hg(0)}{n^2} \sum_{t=1}^n \sqrt{\mathrm{Var}_n X_{t-1}}. \qquad (5)$$

From (2) we have for $\theta < 1$, $\mathrm{Var}_\theta X_t \leq (\sigma_G^2 + \mu_G)(1 - \theta^{2t})(1 - \theta^2)^{-1}$. And for $1 \leq t \leq n$ we have $0 \leq 1 - (1 - h/n^2)^{2t} \leq n^{-1} \exp(2h)$. Now we easily obtain

$$\frac{1}{n^2} \sum_{t=1}^n \sqrt{\mathrm{Var}_n X_{t-1}} \leq \sqrt{\sigma_G^2 + \mu_G} \sqrt{\frac{1}{n^2} \frac{1}{1-(1-h/n^2)^2}} \frac{1}{n} n \sqrt{\frac{\exp(2h)}{n}} \to 0 \qquad \text{as } n \to \infty.$$

A combination with (5) yields (iii). This concludes the proof. $\square$

**Proof of Lemma 3.1.** We obtain, for $h > 0, h_0 \geq 0$, using the inequality $|\log((a+b)/(c+d)) - \log(a/c)| \leq b/a + d/c$ for $a, c > 0$, $b, d \geq 0$, the bound

$$|\mathcal{L}_n(h, h_0) - \tilde{\mathcal{L}}_n(h, h_0)| \leq \sum_{t=1}^n \frac{R_n(X_{t-1}, X_t, h)}{L_n(X_{t-1}, X_t, h)} + \sum_{t=1}^n \frac{R_n(X_{t-1}, X_t, h_0)}{L_n(X_{t-1}, X_t, h_0)} \qquad \mathbb{P}_{1-h_0/n^2}\text{-a.s.} \qquad (6)$$

It is easy to see, since $b_{n,0}(k) = 0$ for $k > 0$ and $g(i) = 0$ for $i < 0$, that for $h_0 > 0$, $\mathcal{L}_n(0, h_0)$ and $\tilde{\mathcal{L}}_n(0, h_0)$ both contain $\log(0)$ exactly $\sum_{t=1}^n 1\{\Delta X_t < 0\}$ times. Also for $\sum_{t=1}^n 1\{\Delta X_t < 0\} = 0$ we have

$$|\mathcal{L}_n(0, h_0) - \tilde{\mathcal{L}}_n(0, h_0)| \leq \sum_{t=1}^n \frac{R_n(X_{t-1}, X_t, h_0)}{L_n(X_{t-1}, X_t, h_0)} \qquad \mathbb{P}_{1-h_0/n^2}\text{-a.s.}$$

Thus if we show that

$$\sum_{t=1}^n \frac{R_n(X_{t-1}, X_t, h')}{L_n(X_{t-1}, X_t, h')} \xrightarrow{p} 0, \qquad \text{under } \mathbb{P}_{1-h_0/n^2}$$

holds for $h' = h$ and $h' = h_0$, the lemma is proved (exclude the case $h' = 0$ and $h_0 > 0$, which need not be considered). We split the expression in the previous display into four non-negative parts (empty sums are by definition equal to 0)

$$\sum_{t=1}^n \frac{R_n(X_{t-1}, X_t, h')}{L_n(X_{t-1}, X_t, h')} = \sum_{t:\Delta X_t \leq -2} \frac{R_n(X_{t-1}, X_t, h')}{L_n(X_{t-1}, X_t, h')} + \sum_{t:\Delta X_t = -1} \frac{R_n(X_{t-1}, X_t, h')}{L_n(X_{t-1}, X_t, h')}$$
$$+ \sum_{t:0 \leq \Delta X_t \leq M} \frac{R_n(X_{t-1}, X_t, h')}{L_n(X_{t-1}, X_t, h')} + \sum_{t:\Delta X_t > M} \frac{R_n(X_{t-1}, X_t, h')}{L_n(X_{t-1}, X_t, h')}.$$

Since $\Delta X_t \leq -2$ implies $X_{t-1} - \vartheta \circ X_{t-1} \geq 2$, (7) implies

$$\sum_{t:\Delta X_t \leq -2} \frac{R_n(X_{t-1}, X_t, h')}{L_n(X_{t-1}, X_t, h')} \xrightarrow{p} 0, \qquad \text{under } \mathbb{P}_{1-h_0/n^2}.$$

Next we treat the terms for which $\Delta X_t = -1$. If $h_0 = 0$, we do not have such terms (under $\mathbb{P}_{1-h_0/n^2}$) and remember that the case $h' = 0$ and $h_0 > 0$ need not be considered. So we only need to consider this term for $h', h_0 > 0$. On the event $\mathcal{A}_n^{h'}$ (see (5) for the definition of this event), an application of (21) yields,

$$\sum_{t:\Delta X_t = -1} \frac{R_n(X_{t-1}, X_t, h')}{L_n(X_{t-1}, X_t, h')} \leq \sum_{t:\Delta X_t = -1} \frac{\sum_{k=3}^{X_{t-1}} b_{X_{t-1}, h'/n^2}(k)}{g(0) b_{X_{t-1}, h'/n^2}(1)}$$
$$\leq 2 \sum_{t=1}^n \frac{(X_{t-1}^3/3!)(h'^3/n^6)(1 - h'/n^2)^{X_{t-1}-3}}{g(0) X_{t-1}(h'/n^2)(1 - h'/n^2)^{X_{t-1}-1}} 1\{X_{t-1} \geq 1\}$$
$$\leq \frac{4h'^2}{3g(0)n^4} \sum_{t=1}^n X_{t-1}^2,$$

since $(1 - h'/n^2)^{-2} \leq 4$ by definition of $\mathcal{A}_n^{h'}$. From (4) and (6) we now obtain

$$\sum_{t:\Delta X_t=-1} \frac{R_n(X_{t-1}, X_t, h')}{L_n(X_{t-1}, X_t, h')} \xrightarrow{p} 0, \qquad \text{under } \mathbb{P}_{1-h_0/n^2}.$$

Next, we analyze the terms for which $0 \leq \Delta X_t \leq M$. We have, by (21), on the event $\mathcal{A}_n^{h'}$,

$$\sum_{t:0\leq\Delta X_t\leq M} \frac{R_n(X_{t-1}, X_t, h')}{L_n(X_{t-1}, X_t, h')} \leq \sum_{t:0\leq\Delta X_t\leq M} \frac{\sum_{k=2}^{X_{t-1}} b_{X_{t-1},h'/n^2}(k) g(\Delta X_t + k)}{g(\Delta X_t) b_{X_{t-1},h'/n^2}(0)}$$

$$\leq \frac{2}{m^*} \sum_{t:0\leq\Delta X_t\leq M} \frac{b_{X_{t-1},h'/n^2}(2)}{b_{X_{t-1},h'/n^2}(0)}$$

$$\leq \frac{4h'^2}{m^* n^4} \sum_{t=1}^n X_{t-1}^2,$$

where $m^* = \min\{g(k) \mid 0 \leq k \leq M\} > 0$. Now (4) and (6) yield the desired convergence,

$$\sum_{t:0\leq\Delta X_t\leq M} \frac{R_n(X_{t-1}, X_t, h')}{L_n(X_{t-1}, X_t, h')} \xrightarrow{p} 0, \qquad \text{under } \mathbb{P}_{1-h_0/n^2}.$$

Finally, we discuss the terms for which $\Delta X_t > M$. If the support of $G$ equals $\{0, \ldots, M\}$, there are no such terms. So we only need to consider the case where the support of $G$ is $\mathbb{Z}_+$. Since $g$ is non-increasing on $\{M, M+1, \ldots\}$, we have, by (21), $R_n(X_{t-1}, X_t, h') \leq 2g(\Delta X_t) b_{X_{t-1},h'/n^2}(2)$ for $X_{t-1} \in A_n^{h'}$, which yields,

$$0 \leq \frac{R_n(X_{t-1}, X_t, h')}{L_n(X_{t-1}, X_t, h')} \leq \frac{2g(\Delta X_t)(X_{t-1}^2/2)h'^2/n^4(1-h'/n^2)^{X_{t-1}-2}}{g(\Delta X_t)(1-h'/n^2)^{X_{t-1}}}$$

$$\leq \frac{4h'^2}{n^4} X_{t-1}^2, \qquad X_{t-1} \in A_n^{h'}.$$

From (4) and (6) it now easily follows that we have

$$\sum_{t:\Delta X_t \geq M} \frac{R_n(X_{t-1}, X_t, h')}{L_n(X_{t-1}, X_t, h')} \xrightarrow{p} 0, \qquad \text{under } \mathbb{P}_{1-h_0/n^2}.$$

This concludes the proof of the lemma. $\square$

**Proof of Lemma 3.2.** We write $S_n^+(h, h_0) = \sum_{t:\Delta X_t \geq 0} \log[1 + U_{tn}^+]$, where

$$U_{tn}^+ = \left( g(\Delta X_t) \left[ \frac{h^2 - h_0^2}{n^4} - 2\frac{h - h_0}{n^2} \right] + X_{t-1} g(\Delta X_t + 1) \left[ \frac{h - h_0}{n^2} - \frac{h^2 - h_0^2}{n^4} \right] \right)$$

$$\times \left( g(\Delta X_t) \left( 1 - \frac{h_0}{n^2} \right)^2 + X_{t-1} g(\Delta X_t + 1) \frac{h_0}{n^2} \left( 1 - \frac{h_0}{n^2} \right) \right)^{-1}.$$

Notice that, for $n$ large enough,

$$U_{tn}^{+2} \leq \left(2\left(g^2(\Delta X_t)\left[\frac{h^2-h_0^2}{n^4} - 2\frac{h-h_0}{n^2}\right]^2 + X_{t-1}^2 g^2(\Delta X_t + 1)\left[\frac{h-h_0}{n^2} - \frac{h^2-h_0^2}{n^4}\right]^2\right)\right)$$
$$\times \left(g^2(\Delta X_t)\left(1 - \frac{h_0}{n^2}\right)^4\right)^{-1}$$
$$\leq \frac{C}{n^4}(X_{t-1}^2 + 1),$$

for some constant $C$, where we used that $e \mapsto g(e+1)/g(e)$ is bounded. From (4) we obtain

$$\lim_{n\to\infty} \mathbb{E}_{1-h_0/n^2} \sum_{t:\Delta X_t \geq 0} U_{tn}^{+2} \leq 0 + \lim_{n\to\infty} \mathbb{E}_{1-h_0/n^2} \frac{C}{n^4} \sum_{t=1}^n X_{t-1}^2 = 0.$$

Hence

$$\sum_{t:\Delta X_t \geq 0} U_{tn}^{+2} \xrightarrow{P} 0, \qquad \text{under } \mathbb{P}_{1-h_0/n^2} \quad \text{and}$$

$$\lim_{n\to\infty} \mathbb{P}_{1-h_0/n^2}\left\{\max_{t:\Delta X_t \geq 0} |U_{tn}^+| \leq 1/2\right\} = 1. \tag{7}$$

Using $\log(1+x) = x + r(x)$, where $r$ satisfies $|r(x)| \leq 2x^2$ for $|x| \leq 1/2$, we obtain from (7),

$$S_n^+(h, h_0) = \sum_{t:\Delta X_t \geq 0} \log[1 + U_{tn}^+] = \sum_{t:\Delta X_t \geq 0} U_{tn}^+ + \mathrm{o}(\mathbb{P}_{1-h_0/n^2}; 1).$$

Thus the problem reduces to determining the asymptotic behavior of $\sum_{t:\Delta X_t \geq 0} U_{tn}^+$. Note that,

$$\sum_{t:\Delta X_t \geq 0} U_{tn}^+ = \sum_{t:\Delta X_t \geq 0} \frac{X_{t-1} g(\Delta X_t + 1)[(h-h_0)/n^2 - (h^2-h_0^2)/n^4]}{g(\Delta X_t)(1-h_0/n^2)^2 + X_{t-1} g(\Delta X_t + 1)(h_0/n^2)(1-h_0/n^2)}$$
$$+ \mathrm{o}(\mathbb{P}_{1-h_0/n^2}; 1).$$

Using that $e \mapsto g(e+1)/g(e)$ is bounded and (4), we obtain

$$\sum_{t:\Delta X_t \geq 0} \left|\frac{X_{t-1} g(\Delta X_t + 1)[(h-h_0)/n^2 - (h^2-h_0^2)/n^4]}{g(\Delta X_t)(1-h_0/n^2)^2 + X_{t-1} g(\Delta X_t + 1)(h_0/n^2)(1-h_0/n^2)}\right.$$
$$\left. - \frac{(h-h_0)}{n^2}\frac{X_{t-1} g(\Delta X_t + 1)}{g(\Delta X_t)}\right|$$
$$\leq \frac{C}{n^4}\sum_{t=1}^n X_{t-1}^2 \xrightarrow{P} 0, \qquad \text{under } \mathbb{P}_{1-h_0/n^2}.$$

Thus the previous three displays and (7) yield

$$S_n^+(h,h_0) = \frac{h-h_0}{n^2}\sum_{t=1}^n X_{t-1}\frac{g(\Delta X_t+1)}{g(\Delta X_t)}1\{\Delta X_t \geq 0, X_{t-1}-\vartheta\circ X_{t-1}\leq 1\}+\mathrm{o}(\mathbb{P}_{1-h_0/n^2};1).$$

Finally, we will show that

$$\frac{1}{n^2}\sum_{t=1}^n X_{t-1}\frac{g(\Delta X_t+1)}{g(\Delta X_t)}1\{\Delta X_t \geq 0, X_{t-1}-\vartheta\circ X_{t-1}\leq 1\}$$

$$\xrightarrow{p} \frac{(1-g(0))\mu_G}{2}, \qquad \text{under } \mathbb{P}_{1-h_0/n^2}, \tag{8}$$

which will conclude the proof. For notational convenience we introduce

$$Z_t = \frac{g(\Delta X_t+1)}{g(\Delta X_t)}1\{\Delta X_t \geq 0, X_{t-1}-\vartheta\circ X_{t-1}\leq 1\}$$

$$= \frac{g(\varepsilon_t+1)}{g(\varepsilon_t)}1\{X_{t-1}-\vartheta\circ X_{t-1}=0\} + \frac{g(\varepsilon_t)}{g(\varepsilon_t-1)}1\{\varepsilon_t \geq 1, X_{t-1}-\vartheta\circ X_{t-1}=1\}.$$

Using that $\varepsilon_t$ is independent of $X_{t-1}-\vartheta\circ X_{t-1}$, we obtain

$$\mathbb{E}_{1-h_0/n^2}[Z_t|X_{t-1}-\vartheta\circ X_{t-1}] = (1-g(0))1\{X_{t-1}-\vartheta\circ X_{t-1}=0\}$$

$$+ 1\{X_{t-1}-\vartheta\circ X_{t-1}=1\}\mathbb{E}\frac{g(\varepsilon_t)}{g(\varepsilon_t-1)}1\{\varepsilon_t \geq 1\},$$

where we used that $\mathbb{E}g(\varepsilon_1+1)/g(\varepsilon_1) = 1-g(0)$ and $\mathbb{E}1\{\varepsilon_1 \geq 1\}g(\varepsilon_1)/g(\varepsilon_1-1) < \infty$, since we assumed that $g$ is eventually decreasing. So we have

$$Z_t - \mathbb{E}_{1-h_0/n^2}[Z_t|X_{t-1}-\vartheta\circ X_{t-1}]$$

$$= \left(\frac{g(\varepsilon_t+1)}{g(\varepsilon_t)} - \mathbb{E}\frac{g(\varepsilon_t+1)}{g(\varepsilon_t)}\right)1\{X_{t-1}-\vartheta\circ X_{t-1}=0\}$$

$$+ \left(\frac{g(\varepsilon_t)}{g(\varepsilon_t-1)}1\{\varepsilon_t \geq 1\} - \mathbb{E}\frac{g(\varepsilon_t)}{g(\varepsilon_t-1)}1\{\varepsilon_t \geq 1\}\right)1\{X_{t-1}-\vartheta\circ X_{t-1}=1\}.$$

From this it is not hard to see that we have

$$\mathbb{E}_{1-h_0/n^2}X_{t-1}(Z_t - \mathbb{E}_{1-h_0/n^2}[Z_t|X_{t-1}-\vartheta\circ X_{t-1}]) = 0,$$

$$\mathbb{E}_{1-h_0/n^2}X_{t-1}(Z_t - \mathbb{E}_{1-h_0/n^2}[Z_t|X_{t-1}-\vartheta\circ X_{t-1}])$$

$$\times X_{s-1}(Z_s - \mathbb{E}_{1-h_0/n^2}[Z_s|X_{s-1}-\vartheta\circ X_{s-1}]) = 0 \qquad \text{for } s<t,$$

$$\mathbb{E}_{1-h_0/n^2}(Z_t - \mathbb{E}_{1-h_0/n^2}[Z_t|X_{t-1}-\vartheta\circ X_{t-1}])^2 \leq C, \tag{9}$$

for $C = 2(\mathrm{Var}(g(\varepsilon_1+1)/g(\varepsilon_1)) + \mathrm{Var}(1_{\{\varepsilon_t \geq 1\}} g(\varepsilon_1)/g(\varepsilon_1-1)))$, which is finite by Assumption 3.1. Thus, by (4), it follows that

$$\mathbb{E}_{1-h_0/n^2}\left(\frac{1}{n^2}\sum_{t=1}^n X_{t-1}(Z_t - \mathbb{E}_{1-h_0/n^2}[Z_t|X_{t-1} - \vartheta \circ X_{t-1}])\right)^2$$

$$= \frac{1}{n^4}\sum_{t=1}^n \mathbb{E}_{1-h_0/n^2} X_{t-1}^2 (Z_t - \mathbb{E}_{1-h_0/n^2}[Z_t|X_{t-1} - \vartheta \circ X_{t-1}])^2$$

$$\leq \frac{C}{n^4}\sum_{t=1}^n \mathbb{E}_{1-h_0/n^2} X_{t-1}^2 \to 0.$$

Hence (8) is equivalent to

$$\frac{1}{n^2}\sum_{t=1}^n X_{t-1}\mathbb{E}_{1-h_0/n^2}[Z_t|X_{t-1} - \vartheta \circ X_{t-1}] \xrightarrow{p} \frac{(1-g(0))\mu_G}{2}, \qquad \text{under } \mathbb{P}_{1-h_0/n^2}. \quad (10)$$

Since, by (4),

$$\frac{1}{n^2}\sum_{t=1}^n \mathbb{E}_{1-h_0/n^2} X_{t-1} 1\{X_{t-1} - \vartheta \circ X_{t-1} = 1\} = \frac{h_0}{n^4}\sum_{t=1}^n \mathbb{E}_{1-h_0/n^2} X_{t-1}^2 \left(1 - \frac{h_0}{n^2}\right)^{X_{t-1}-1}$$

$$\leq \frac{h_0}{n^4}\sum_{t=1}^n \mathbb{E}_{1-h_0/n^2} X_{t-1}^2 \to 0,$$

we have, using (7),

$$\left|\frac{1}{n^2}\sum_{t=1}^n X_{t-1}\mathbb{E}_{1-h_0/n^2}[Z_t|X_{t-1} - \vartheta \circ X_{t-1}] - \frac{1-g(0)}{n^2}\sum_{t=1}^n X_{t-1}\right|$$

$$\leq \left|\mathbb{E}\frac{g(\varepsilon_t)}{g(\varepsilon_t-1)}1\{\varepsilon_t \geq 1\} - (1-g(0))\right|\frac{1}{n^2}\sum_{t=1}^n X_{t-1}1\{X_{t-1} - \vartheta \circ X_{t-1} = 1\}$$

$$+ \frac{1-g(0)}{n^2}\sum_{t=1}^n X_{t-1}1\{X_{t-1} - \vartheta \circ X_{t-1} \geq 2\} \xrightarrow{p} 0, \qquad \text{under } \mathbb{P}_{1-h_0/n^2}.$$

We conclude (10), which finally concludes the proof of the lemma. □

**Proof of Lemma 3.3.** First we consider $h = 0$. From the definition of $S_n^-(0, h_0)$ we see that $S_n^-(0, h_0) = 0$ if $\sum_{t=1}^n 1\{\Delta X_t < 0\} = 0$ (since an empty sum equals zero by definition). And if $\sum_{t=1}^n 1\{\Delta X_t < 0\} \geq 1$, we have $S_n^-(0, h_0) = -\infty$ (since $W_{tn}^- = -\infty$ for $h = 0$). This concludes the proof for $h = 0$.

So we now consider $h > 0$. We rewrite

$$W_{tn}^- = \log\left[\left(\frac{h}{h_0}\left(\frac{1-h/n^2}{1-h_0/n^2}\right) + \frac{X_{t-1}-1}{2n^2}\frac{h^2 g(1)}{g(0)h_0(1-h_0/n^2)}\right) \right.$$
$$\left. \times \left(1 + \frac{X_{t-1}-1}{2n^2}\frac{h_0 g(1)}{g(0)(1-h_0/n^2)}\right)^{-1}\right].$$

By (7), the proof is finished if we show that

$$\sum_{t:\Delta X_t=-1}\left|W_{tn}^- - \log\left[\frac{h}{h_0}\right]\right| \xrightarrow{p} 0, \qquad \text{under } \mathbb{P}_{1-h_0/n^2}.$$

Using the inequality $|\log((a+b)/(c+d)) - \log(a/c)| \le b/a + d/c$ for $a, c > 0$, $b, d \ge 0$, we obtain

$$\left|W_{tn}^- - \log\left[\frac{h}{h_0}\right]\right|$$
$$\le \left|W_{tn}^- - \log\left[\frac{h}{h_0}\left(\frac{1-h/n^2}{1-h_0/n^2}\right)\right]\right| + \mathrm{O}(n^{-2})$$
$$\le \frac{X_{t-1}-1}{2n^2}\left[\frac{h^2 g(1)}{g(0)h_0(1-h_0/n^2)}\left(\frac{h}{h_0}\left(\frac{1-h/n^2}{1-h_0/n^2}\right)\right)^{-1} + \frac{h_0 g(1)}{g(0)(1-h_0/n^2)}\right] + \mathrm{O}(n^{-2}).$$

Hence, it suffices to show that

$$\sum_{t:\Delta X_t=-1}\frac{X_{t-1}}{n^2} \xrightarrow{p} 0, \qquad \text{under } \mathbb{P}_{1-h_0/n^2}.$$

Note first that we have, by (7),

$$0 \le \frac{1}{n^2}\sum_{t=1}^n X_{t-1}\mathbf{1}\{\Delta X_t = -1\} = \frac{1}{n^2}\sum_{t=1}^n X_{t-1}\mathbf{1}\{\Delta X_t = -1, \varepsilon_t = 0\} + \mathrm{o}(\mathbb{P}_{1-h_0/n^2}; 1).$$

We show that the expectation of the first term on the right-hand side in the previous display converges to zero, which will conclude the proof. We have, by (4),

$$\lim_{n\to\infty}\frac{1}{n^2}\sum_{t=1}^n \mathbb{E}_{1-h_0/n^2}X_{t-1}\mathbf{1}\{\Delta X_t = -1, \varepsilon_t = 0\}$$
$$= \lim_{n\to\infty}\frac{h_0}{n^4}\sum_{t=1}^n \mathbb{E}_{1-h_0/n^2}g(0)X_{t-1}^2\left(1-\frac{h_0}{n^2}\right)^{X_{t-1}-1}$$
$$\le \lim_{n\to\infty}\frac{h_0 g(0)}{n^4}\sum_{t=1}^n \mathbb{E}_{1-h_0/n^2}X_{t-1}^2 = 0,$$

which concludes the proof of the lemma. □

# Acknowledgements

The authors thank Marc Hallin and Johan Segers for useful discussions and suggestions. Comments by the conference participants at the Econometric Society EM 2006, Prague Stochastics 2006 and Faro $EC^2$ as well as seminar participants at Universiteit van Amsterdam and Université Libre de Bruxelles are kindly acknowledged. Furthermore, the authors are grateful to the referees for very helpful comments and suggestions.